\documentclass[12pt,leqno,fleqn,epsfig]{article}
\usepackage{amssymb,epsfig, amsmath, amsthm}
\usepackage{mathrsfs}       

\newcommand{\R}{\mathbb{R}}

\newcommand{\N}{\mathbb{N}}

\newcommand{\ba}{\boldsymbol a}
\newcommand{\bb}{\boldsymbol b}
\newcommand{\bc}{\boldsymbol c}
\newcommand{\bd}{\boldsymbol d}
\newcommand{\bfe}{\boldsymbol e}
\newcommand{\bbf}{\boldsymbol f}
\newcommand{\bg}{\boldsymbol g}
\newcommand{\bh}{\boldsymbol h}
\newcommand{\bi}{\boldsymbol i}
\newcommand{\bj}{\boldsymbol j}
\newcommand{\bk}{\boldsymbol k}
\newcommand{\bl}{\boldsymbol l}
\newcommand{\bm}{\boldsymbol m}
\newcommand{\bn}{\boldsymbol n}
\newcommand{\bo}{\boldsymbol o}
\newcommand{\bp}{\boldsymbol p}
\newcommand{\bq}{\boldsymbol q}
\newcommand{\br}{\boldsymbol r}
\newcommand{\bs}{\boldsymbol s}
\newcommand{\bt}{\boldsymbol t}
\newcommand{\bu}{\boldsymbol u}
\newcommand{\bv}{\boldsymbol v}
\newcommand{\bw}{\boldsymbol w} 
\newcommand{\bx}{\boldsymbol x}
\newcommand{\by}{\boldsymbol y}
\newcommand{\bz}{\boldsymbol z}
\newcommand{\buu}{ \underline{\boldsymbol u}}

\newcommand{\bA}{\boldsymbol A}
\newcommand{\bB}{\boldsymbol B}
\newcommand{\bC}{\boldsymbol C}
\newcommand{\bD}{\boldsymbol D}
\newcommand{\bE}{\boldsymbol E}
\newcommand{\bF}{\boldsymbol F}
\newcommand{\bG}{\boldsymbol G}
\newcommand{\bH}{\boldsymbol H}
\newcommand{\bI}{\boldsymbol I}
\newcommand{\bJ}{\boldsymbol J}
\newcommand{\bK}{\boldsymbol K}
\newcommand{\bL}{\boldsymbol L}
\newcommand{\bM}{\boldsymbol M}
\newcommand{\bO}{\boldsymbol O}
\newcommand{\bP}{\boldsymbol P}
\newcommand{\bQ}{\boldsymbol Q}
\newcommand{\bR}{\boldsymbol R}
\newcommand{\bS}{\boldsymbol S}
\newcommand{\bT}{\boldsymbol T}
\newcommand{\bU}{\boldsymbol U}
\newcommand{\bV}{\boldsymbol V}
\newcommand{\bW}{\boldsymbol W}
\newcommand{\bX}{\boldsymbol X}
\newcommand{\bY}{\boldsymbol Y}
\newcommand{\bZ}{\boldsymbol Z}

\newcommand{\bfvarphi}{\boldsymbol\varphi}

\newcommand{\bPhi}{\boldsymbol \Phi}
\newcommand{\bfomega}{\boldsymbol \omega}

\newcommand{\bLambda}{\boldsymbol \Lambda}

\newcommand{\mes}{\operatorname{\rm mes}}

\newcommand{\osc}{\operatorname*{osc}}

\newcommand{\supp}{\operatorname*{supp}}

\newcommand{\dist}{\operatorname*{dist}}

\newcommand{\const}{\operatorname*{const}}

\newcommand{\be}{\begin{equation}}
\newcommand{\ee}{\end{equation}}
\newcommand{\bea}{\begin{eqnarray}}
\newcommand{\eea}{\end{eqnarray}}
\newcommand{\bean}{\begin{eqnarray*}}
\newcommand{\eean}{\end{eqnarray*}}

\newcommand{\var}{\varepsilon}

\newcommand{\intl}{\int\limits}
\newcommand{\suml}{\sum\limits}

\newcommand{\Beweisende}{\rule{0.2cm}{0.2cm}}

\newcommand{\D}{\displaystyle}

\newcommand{\intmw}{{\int\hspace{-830000sp}-\!\!}}

\newcounter{secnum}

\newtheorem{thm}{Theorem}[section]
\newtheorem{cor}[thm]{Corollary}
\newtheorem{lem}[thm]{Lemma}

\theoremstyle{definition}
\newtheorem{defin}[thm]{Definition}
\newtheorem{rem}[thm]{Remark}

\title{On partial regularity for the steady Hall magnetohydrodynamics system}
 \author{Dongho Chae$^*$  and J\"{o}rg Wolf $^\dagger$\\
\ \\
 $*$Department of Mathematics\\
Chung-Ang University\\
 Seoul 156-756, Republic of Korea\\
 e-mail: dchae@cau.ac.kr\\
and \\
$\dagger$Department of Financial Engineering\\
 Ajou University\\                                           
Suwon 443-749, Republic of Korea\\
e-mail: jwolf@math.hu-berlin.de}
\date{}
\begin{document}
\maketitle
\begin{abstract}
We study partial regularity of suitable weak solutions of the steady Hall magnetohydrodynamics equations in a domain $\Omega \subset \Bbb R^3$.  In particular we prove that the set of possible singularities of the suitable weak solution has Hausdorff dimension at most one.  Moreover, in the case $\Omega=\Bbb R^3$, we show that the 
set of possible singularities is compact.
\\
\ \\
\noindent{\bf AMS Subject Classification Number:} 35Q35, 35Q85,76W05\\
  \noindent{\bf
keywords:} steady Hall-MHD equations, parial regularity

\end{abstract}

\section{Introduction}
\label{sec:-1}
\setcounter{secnum}{\value{section} \setcounter{equation}{0}
\renewcommand{\theequation}{\mbox{\arabic{secnum}.\arabic{equation}}}}
The resistive incompressible Hall magnetohydrodynamics(Hall-MHD) is  described by he following equations:
$$
\left\{ \aligned 
&\frac{\partial\bu}{\partial t} +\bu\cdot\nabla \bu+\nabla p=(\nabla\times \bB)\times \bB +\nu\Delta \bu +\bbf,\\
&\frac{\partial\bB}{\partial t}-\nabla\times (\bu\times \bB)+\nabla\times ((\nabla\times \bB)\times \bB)
=\mu \Delta B +\nabla \times \bg,\\
&\nabla \cdot \bu=0,\,\,\nabla \cdot \bB=0,\\ \endaligned \right.
$$
where  3D vector fields $\bu=\bu(x,t),  \bB=\bB(x,t)$ are the fluid velocity and the magnetic field respectively.  The scalar field $p=p(x,t)$ is the pressure, while the positive constants $\nu$ and $\mu$ represent the viscosity and the magnetic resistivity respectively. The given vector fields $\bbf$ and $\nabla \times \bg$ are  external forces  on the magnetically charged fluid flows.  Historically, the Hall-MHD system was first considered by Lighthill(\cite{lig}). Compared with the usual MHD system,  the  Hall-MHD system contains the extra term $\nabla\times ((\nabla\times \bB)\times \bB)$, called the Hall term. The inclusion of this  term is essential in understanding the problem of magnetic reconnection, which corresponds to the change of the topology of magnetic field lines.  This phenomena of magnetic reconnection is really observed, for example, in  space plasma(\cite{for, hom}), star formation(\cite{war}) and neutron star(\cite{sha}). 
For the other physical features related to the Hall-MHD we refer \cite{miu, sim},  while for a comprehensive review of the physical aspect of the equations we refer \cite{pol}.  Since the Hall term involves the second order derivative of the magnetic field, it becomes important when the magnetic shear is very large, and this occurs during the reconnection procedure. In the laminar flows this term is small compared with the other term, and can be neglected, which is the case of the usual MHD. 

\hspace{0.5cm}
Since  the Hall term is quadratically nonlinear,   containing the second order derivative, it causes major difficulties in the mathematical study of the Hall-MHD system,  and only recently the rigorous  results on the Cauchy problem appeared. In \cite{ach} the authors 
proved the global existence of weak solutions,  while the local in time well-posedness  as well as the global in time well-posedness for small initial data was proved  in \cite{cha1}. This later result was refined in \cite{cha2}. In the case of $\mu=0$ it is proved in \cite{cha4} that the Cauchy problem is not globally in time well-posed, rigorously verifying the numerical experiment of \cite{dre}. For a special axially symmetric initial data the authors of \cite{fan} proved the global well-posedness of the system. In \cite{cha3} the long time behaviors of the solution were also studied.
Since the Hall-MHD system has more complicated structure than the usual MHD system and the Navier-Stokes equations, the study of full regularity of weak solutions would be extremely difficult. Therefore, it might be reasonable to begin with the partial regularity, similarly to the case of the Navier-Stokes equations, the partial regularity of which was studied e.g.  in \cite{sch, caf, lad, lin, wol}.  In the time-dependent problem, mainly due to the difficulty of defining the correct localized energy inequality for a suitable weak solutions
we concentrate the partial regularity problem of the steady Hall-MHD system. Contrary to the case of the Navier-Stokes equations and the usual MHD system the full regularity of the steady weak solutions 
is difficult to deduce.  Instead,  we prove that the set of possible singularity of the steady suitable weak solution of the Hall-MHD system has Hausdorff dimension at most $1$(see Remark 5.2 and Theorem 5.3 below). Moreover, for a steady suitable weak solution on $\Bbb R^3$ the set of possible singularity is 
a compact set(see Corollary 7.3 below). The partial regularity of the  time dependent problem will be studied elsewhere.

%
%
\section{Weak solution  and higher regularity of $\bu$ } 
\setcounter{secnum}{\value {section}
\setcounter{equation}{0}
\renewcommand{\theequation}{\mbox{\arabic{secnum}.\arabic{equation}}}}

We consider the following steady Hall-MHD system in $\R^3$.
\begin{align}
 (\bu\cdot \nabla )\bu - \Delta \bu &= -\nabla p +
 (\nabla \times \bB) \times  \bB + \bbf,
\label{1.2}
\\
\nabla \times (\bB\times \bu) - \Delta \bB &= 
-\nabla \times ((\nabla \times \bB) \times  \bB)+  \nabla \times \bg,
\label{1.3}
\\
\nabla \cdot \bu &=0, \,\, \nabla \cdot \bB= 0.
\label{1.1}  
\end{align}
We note that we set  $\mu=\nu=1$ for convenience.
As $(\bu\cdot \nabla )\bu  = (\nabla \times \bu )\times \bu+ \frac {1} {2} |\bu|^2$, 
\eqref{1.2} turns into 
\be\label{1.4}
 (\nabla \times \bu )\times \bu - \Delta \bu = 
-\nabla \Big(p + \frac {|\bu|^2} {2}\Big)+ (\nabla \times \bB) \times  \bB + \bbf \quad \mbox{in}\quad \R^3.
\ee

Applying $\nabla \times $ to the both sides of the above, we get
\be\label{1.5}
\nabla \times (\bfomega\times \bu) - \Delta \bfomega = 
\nabla \times ((\nabla \times \bB) \times  \bB) + \nabla \times \bbf\quad \mbox{in}\quad \R^3,
\ee
where $\bfomega $  stands for the vorticity $\nabla \times \bu$. 
Taking the sum of \eqref{1.3} and \eqref{1.5}, we are led to 
\be\label{1.6}
\nabla \times (\bV\times \bu) - \Delta \bV =  \nabla \times (\bbf + \bg)  \quad \mbox{in}\quad \R^3,
\ee
where
\be\label{1.7}
\bV = \bB+\bfomega. 
\ee
As $\nabla \cdot \bV=0$, there exists a solenoidal potential $\bv$ such that $\nabla \times \bv= \bV$. 
From \eqref{1.6} we deduce that $\bv$ solves the system in $\R^3$,
\begin{align}
 \nabla \cdot  \bv &=0,
\label{1.8}
\\
 (\bv \cdot \nabla) \bv - \Delta \bv &= 
- \nabla \pi + (\nabla \times \bv) \times \bb + \bbf+ \bg,
\label{1.9}
\end{align}
where $\bb= \bv -\bu$.  Clearly, $\nabla \times \bb = \bB$. 

 \begin{defin}
Let $\bbf\in  L^{6/5} $ and $\bg \in  L^2$.  We say $(\bu, p, \bB) \in  \dot{W}^{ 1,\, 2}\times L^{2}_{\rm loc}\times  \dot{W}^{ 1,\, 2}$  is a 
{\it weak solution} to \eqref{1.1}--\eqref{1.3} if  
\begin{align}
  &\intl_{\R^3} \nabla \bu : \nabla \bfvarphi - \bu \otimes \bu : \nabla \bfvarphi  dx 
\cr
&\quad = \intl_{\R^3} p \nabla \cdot  \bfvarphi  dx  +\intl_{\R^3} ((\nabla \times \bB) \times \bB) \cdot \bfvarphi  dx  + \intl_{\R^3} \bbf\cdot \bfvarphi  dx,   
\label{2.1}
\\
  &\intl_{\R^3} \nabla \bB : \nabla \bfvarphi + \bB \times \bu \cdot \nabla \times \bfvarphi  dx 
\cr
&\quad = -\intl_{\R^3} ((\nabla \times \bB) \times \bB) \cdot \nabla \times \bfvarphi  dx  
+ \intl_{\R^3} \bg \cdot \nabla \times  \bfvarphi  dx  
\label{2.2}
\end{align} 
for all $\bfvarphi \in  C^\infty _{\rm c} $.  Here $\dot{W}^{ 1,\, 2} $ stands for the homogeneous Sobolev space.  

\end{defin}

\begin{rem}
Since  $\bB \in  \dot{W}^{ 1,\, 2}\hookrightarrow L^{6}$ , by the Calder\'on-Zygmund inequality and 
 Sobolev's embedding theorem we infer 
$ \bb \in  L^\infty_{ \rm loc} $.  In particular, having $\bV\in L^2_{\rm loc} $  we easily deduce 
that $\bv \in  W^{ 1,\, 2}_{\rm loc} $ and $\pi \in  L^2_{\rm loc} $.   Furthermore, 
$(\bv, \pi )$  satisfies the following integral identity for all $ \bfvarphi\in C^{\infty}_{\rm c} $
\begin{align}
  &\intl_{\R^3} \nabla \bv : \nabla \bfvarphi - \bv \otimes \bv: \nabla \bfvarphi  dx 
\cr
&\quad = \intl_{\R^3} \pi  \nabla \cdot  \bfvarphi  dx  +\intl_{\R^3} ((\nabla \times \bv) \times \bb) \cdot \bfvarphi  dx  +\intl_{\R^3} (\bbf+\bg)\cdot \bfvarphi  dx.  
\label{2.3}
\end{align}
\end{rem}

\hspace*{0.5cm}
By using standard regularity methods we get the following   

\begin{thm} 
Let $\bbf\in  L^{6/5}$ and $\bg \in L^2$. Let $(\bu, p , \bB)\in  \dot{W}^{ 1,\, 2}\times 
L ^2_{\rm loc}  \times \dot{W}^{ 1,\, 2}$ 
be a weak solution to \eqref{1.1}--\eqref{1.3}.  Suppose, $\bbf, \bg\in  L^q_{\rm loc}(\Omega )$ for 
some $\frac {6} {5}<q<  +\infty  $ and for an open set $\Omega \subset \R^3$.  Then 
\be\label{2.4}
\bV \in  W^{1,\, q}_{\rm loc}(\Omega ),\quad \bv \in  W^{2,\, q}_{\rm loc}(\Omega ),\quad 
\bu \in  W^{2,\, q\wedge 2}_{\rm loc}(\Omega ),
\ee
where $q\wedge 2=\min\{ q,2\}$.
\end{thm}

{\it Proof}  First, assume $\frac {3} {2} \le q < +\infty $.  Since $\bu\in  L^6$ and $\bV\in L^2_{ \rm loc}$ , we see that 
$-\bV \times \bu + \bbf + \bg \in  L^{3/2}_{\rm loc}(\Omega ) $.   Observing \eqref{1.6},  by the aid of 
Calder\'{o}n-Zygmund's inequality and Sobolev's embedding theorem 
we find $\bV \in  W^{ 1,\, 3/2}_{\rm loc}(\Omega ) \subset 
L^3_{\rm loc}(\Omega )$.   As $\bfomega = \bV- \bB$, this implies  
$\bfomega \in  L^3_{\rm loc}(\Omega ) $.  Taking into account $\nabla \cdot \bu=0$, we get 
$\bu \in  W^{ 1,\, 3}_{\rm loc}(\Omega ) $. Once more appealing to Sobolev's embedding therorem,  we 
obtain  
$\bu \in  L^s_{\rm loc}(\Omega ) $ for all $1\le s <+\infty $, and thus 
$-\bV \times \bu + \bbf + \bg \in  L^{q}_{\rm loc} (\Omega )$.  Again applying 
Calder\'{o}n-Zygmund's inequality, we see that  $\bV \in  W^{ 1,\, q}_{\rm loc}(\Omega )$. 
As $\nabla \cdot \bv =0$ from the last statement we infer $\bv \in  W^{ 2,\, q}_{\rm loc}(\Omega )$. 
Finally,  recalling $\bfomega = \bV- \bB$ and $\bB \in  W^{ 1,\, 2}_{\rm loc}(\Omega ) $,  
we   obtain $\bfomega \in  W^{ 1,\, q\wedge 2}_{\rm loc} (\Omega )$.  

\hspace*{0.5cm}
In case $\frac {6} {5} < q < \frac {3} {2}$, we immediately get 
$-\bV \times \bu + \bbf + \bg \in  L^{q}_{\rm loc}(\Omega )$, and the assertion can be proved as in the 
previous case. 
\hfill \Beweisende

\section{Caccioppoli-type inequality  for $\bB$}
\label{sec:-2}
\setcounter{secnum}{\value{section} \setcounter{equation}{0}
\renewcommand{\theequation}{\mbox{\arabic{secnum}.\arabic{equation}}}}

Suppose $\bbf, \bg \in L^2_{\rm loc}(\Omega ) $  for some open set $\Omega \subset \R^3$. 
As it has been proved in Section\,2 we get $\bu \in  W^{ 2,\, 2}_{\rm loc} (\Omega )$ 
if $(\bu, p, \bB)$ is a weak solution to the steady Hall-MHD system. By means of Sobolev's embedding 
theorem this implies $\bu\in L^\infty _{\rm loc}(\Omega )$.  Accordingly, 
$-\bB\times  \bu + \bg \in L^2_{\rm loc}(\Omega ) $.  

\hspace*{0.5cm}
Thus, for the sake of generality in the present and next section 
we  study the local regularity for the following general 
model.  Let $\Omega \subset \R^3$ be a domain.  We consider the system  
\be\label{3.0}
- \Delta \bB = - \nabla \times ((\nabla \times \bB )\times \bB) + 
\nabla \times \bF \quad \mbox{in} \quad \Omega. 
\ee

\hspace{0.5cm}
We start our discussion with the following notion of a weak solution to \eqref{3.0}.

\begin{defin}
Let $\bF\in  L^2_{\rm loc}(\Omega ) $.  (i)  A vector function $\bB \in  W^{ 1,\, 2}_{\rm loc} (\Omega )$ is said to be a {\it weak solution to }\eqref{3.0} if 
\be\label{3.0-1}
\intl_{\Omega} \nabla \bB: \nabla \bfvarphi  dx  = 
- \intl_{\Omega }  (\nabla \times \bB)  \times \bB \cdot \nabla \times \bfvarphi  dx  + 
\intl_{\Omega} \bF \cdot \nabla \times \bfvarphi  dx  
\ee  
for all $\bfvarphi \in  C^\infty _{\rm c}(\Omega ) $. 

\hspace*{0.5cm}
(ii)  A weak solution $\bB$ to \eqref{3.0} is called a {\it suitable weak solution to } \eqref{3.0} 
if, in addition, the following local energy inequality holds:
\begin{align}
 & \intl_{\Omega} \phi  |\nabla \bB|^2 dx 
\cr
&\le  
 \frac {1} {2} \intl_{\Omega} \Delta \phi |\bB|^2  dx     
- \intl_{\Omega} (\nabla \times \bB) \times  \bB \cdot   (\bB \times  \nabla \phi )       dx  
\cr
&\qquad\qquad+ \intl_{\Omega}( \phi \bF \cdot  \nabla \times \bB + \bF \cdot \bB \times \nabla \phi  ) dx  
\label{3.0-3}
\end{align}
for all non-negative $\phi \in C^\infty _{\rm c}(\Omega ) $.

\end{defin}
\begin{rem}
Let $\bB$ be a suitable weak solution to \eqref{3.0}. Then  for every constant vector 
$\bLambda \in \R^3$  there holds 
\begin{align}
 & \intl_{\Omega} \phi  |\nabla \bB|^2 dx 
\cr
&\le  
 \frac {1} {2} \intl_{\Omega} \Delta \phi |\bB- \bLambda |^2  dx     
- \intl_{\Omega} (\nabla \times \bB) \times  \bB \cdot 
  ((\bB - \bLambda) \times  \nabla \phi )       dx  
\cr
&\qquad\qquad+ \intl_{\Omega} \left\{\phi \bF \cdot  \nabla \times \bB + 
\bF \cdot (\bB- \bLambda ) \times \nabla \phi  \right\} dx  
\label{3.0-4}
\end{align}
for all non-negative $\phi \in C^\infty _{\rm c}(\Omega ) $. 
This can be readily seen by combining \eqref{3.0-3} and \eqref{3.0-1} with $\bfvarphi = \phi \bLambda $.
\end{rem}

\hspace*{0.5cm}
Now, we state the following Caccioppoli-type inequality 

\begin{lem}
Let $\bF\in  L^2_{\rm loc} (\Omega ) $. Let $\bB \in  W^{ 1,\, 2}_{\rm loc} (\Omega )$ be a suitable weak solution to \eqref{3.0}. Then for every ball $B_r= B_r(x_0)\subset \subset  \Omega $ 
and $0< \rho < r$ there holds 
\begin{align}
& \frac {1} {r} \intl_{B_{\rho  } }|\nabla  \bB|^2  dx  
\cr
&\quad\le \frac {cr^2} {(r-\rho )^2} (1+ |\bB_{r, x_0} |^2) \intmw_{B_r} |\bB- \bB_{r, x_0} |^2  dx  
+ \frac {cr^2} {(r-\rho )^2} \intmw_{B_r} |\bB- \bB_{r,x_0} |^4  dx
\cr 
&\qquad\qquad+ \frac {c} {r} \|\bF\|_{2, B_r}^2,    
\label{3.5}
\end{align}
where $\bB_{r, x_0}$ stands for the mean value 
\[
\intmw_{B_r}\bB dx := \frac{1}{{\rm meas}(B_r)} \intl_{B_r(x_0)} \bB  dx, 
\]
 and $c=\const>0$ denotes a universal constant.

\end{lem}

{\it Proof}  Let $B_r = B_r(x_0)\subset \subset  \Omega $  be a fixed ball.  
Given $\rho \in (0, r)$, we consider a cut-off function $\zeta$ defined as 
 $\zeta \in  C^\infty _{\rm c}(B_{r} )$ such that $0\le \zeta  \le 1$  in $\R^3$, $\zeta \equiv 1$ on $B_{\rho } $ and 
$|\nabla \zeta |^2 + |\nabla ^2 \zeta |\le c(r-\rho )^{-2} $ in $\R^3$.
From \eqref{3.0-4} with $\phi = \zeta ^2$ and $\bLambda = \bB_{r, x_0} $   we obtain 
the following Caccioppoli-type inequality
\begin{align}
& \intl_{B_r}  \zeta ^2|\nabla  \bB|^2  dx  
\cr
&\quad\le \frac {c} {(r-\rho )^2} \intl_{B_r} |\bB- \bB_{r, x_0} |^2  dx  +  c\intl_{B_r} |\bF|^2  dx
\cr
&\qquad\qquad + \frac {c} {r-\rho }\intl_{B_r} \zeta |\nabla \bB|\, |\bB|\,|\bB- \bB_{r, x_0}|dx
\cr
&  \quad =  \frac {c} {(r-\rho )^2} \intl_{B_r} |\bB- \bB_{r, x_0} |^2  dx   +   
c\intl_{B_r} |\bF|^2  dx  +J. 
\label{3.1}
\end{align}
Applying H\"older's and Young's inequality, we estimate
\begin{align*}
J& \le \frac {c} {(r-\rho )^2}\intl_{B_r} |\bB|^2|\bB- \bB_{r, x_0}|^2 dx 
+\frac {1} {2} \intl_{B_r}  \zeta ^2|\nabla  \bB|^2  dx  
\\     
& \le \frac {c} {(r-\rho )^2} |\bB_{r, x_0}   |^2 \intl_{B_r} |\bB- \bB_{r, x_0} |^2 dx 
+ \frac {c} {(r-\rho )^2}\intl_{B_r} |\bB- \bB_{r, x_0} |^4 dx
\\
&\qquad\qquad+\frac {1} {2} \intl_{B_r}  \zeta ^2|\nabla  \bB|^2  dx.
\end{align*}
Inserting the above estimate of $J$ into the right-hand side of \eqref{3.1}, and dividing the resulting  
estimate by $r$,  we are led to 
\begin{align}
& \frac {1} {r} \intl_{B_{\rho  } }|\nabla  \bB|^2  dx  
\cr
&\quad\le \frac {cr^2} {(r-\rho )^2} (1+ |\bB_{r, x_0} |^2) \intmw_{B_r} |\bB- \bB_{r, x_0}|^2  dx  
+ \frac {cr^2} {(r-\rho )^2} \intmw_{B_r} |\bB- \bB_{r, x_0}|^4  dx
\cr 
&\qquad\qquad+ \frac {c} {r}\|\bF\|^2_{2, B_r}.    
\label{3.3-1}
\end{align}
Whence, the claim. \hfill \Beweisende

\begin{rem}
Setting 
\begin{align*}
E(r)= E(r, x_0) &=   \bigg(\intmw_{B_r(x_0)} |\bB-  \bB_{r; x_0}  |^4  dx\bigg)^{1/4}\quad 
0<r< \dist(x_0,\partial \Omega ) ,
\end{align*}
\eqref{3.5} becomes
\begin{align}
\bigg(\frac {1} {r} \intl_{B_\rho  }  |\nabla  \bB|^2  dx\bigg)^{1/2}  \le 
\frac {cr} {r-\rho } \Big\{(1+|\bB_{r, x_0} |) E(r) + E(r)^2\Big\} +  
\frac {c} {r^{1/2} }  \|\bF\|_{2, B_r}.
\label{3.4-1}
\end{align}
\end{rem}

%
%
\section{Blow-up} 
\setcounter{secnum}{\value {section}
\setcounter{equation}{0}
\renewcommand{\theequation}{\mbox{\arabic{secnum}.\arabic{equation}}}}

\vspace{0.3cm}
We begin our discussion with the following fundamental estimate for solutions to the model problem, which will be used in the blow-up lemma below.

\vspace{0.3cm}
\begin{lem}
Let $\bLambda \in \R^3$. Let $\bW\in L^4(B_1)\cap  W^{ 1,\, 2}_{\rm loc} (B_1)$ be a 
weak solution to 
\be\label{4.1}
-\Delta \bW = - \nabla \times ((\nabla \times \bW) \times  \bLambda) 
\quad \mbox{in}\quad B_1,
\ee
i.\,e.   
\be\label{4.2}
\intl_{B_1} \nabla \bW : \nabla \bPhi    dx =   - \intl_{B_1} ((\nabla \times \bW) \times \bLambda) \cdot  
\nabla \times \bPhi    dx
\ee
for all $\bPhi \in  W^{ 1,\, 2}(B_1)$ with $\supp (\bPhi ) \subset \subset  B_1$. Then, 
\be\label{4.3}
\bigg(\intmw_{B_\tau } |\bW- \bW_{B_\tau }|^4 dx  \bigg)^{1/4} \le C_0 \tau (1+ |\bLambda  |^{3} ) 
\bigg(\intmw_{B_1} |\bW- \bW_{B_1 }|^4 dx\bigg)^{1/4} \quad \forall\, 0<\tau <1,
\ee
where $C_0>0$ denotes a universal constant. 
\end{lem}

{\it Proof}  Since the assertion is trivial for $\frac {1} {2} < \tau < 1$, we may assume that $0< \tau < \frac {1} {2}$. 
Let $\zeta \in  C^\infty _{\rm c} (B_1)$  be a suitable cut-off function for $B_{1/2}$, i.\,e. $ 0 \le \zeta \le 1$ 
in $ B_1$ and $ \zeta \equiv 1$ on $ B_{ 1/2}$. In \eqref{4.2}  inserting  the admissible test function  $\bPhi =\zeta ^{2m} (\bW-\bW_{B_1}) $\, 
$(m\in \N)$,  by using Cauchy-Schwarz's inequality along with 
Young's inequality we obtain 
\be\label{4.4}
\intl_{B_1} \zeta ^{2m} |\nabla \bW|^2 dx \le 
c (1+ |\bLambda |^2) \intl_{B_1} \zeta ^{2m-2} |\bW- \bW_{B_1} |^2 dx.     
\ee  
If $\bW$ is smooth in $B_1$, since \eqref{4.1} is  a linear system,  the same inequality holds for $D^\alpha \bW $ in place of $\bW$ for  any multi-index $\alpha $. 
By a standard mollifying argument together with Sobolev's embedding theorem 
we see that $\bW$ is smooth. In particular, in \eqref{4.4} putting $m=3$ it follows that 
\be\label{4.5}
\intl_{B_1} \zeta ^{6} |D^\alpha \bW|^2 dx \le 
c (1+ |\bLambda |^6) \intl_{B_1} |\bW- \bW_{B_1} |^2 dx\quad\forall\, |\alpha |\le 3.     
\ee  
By means of Sobolev's embedding theorem and Jensen's inequality we get 
\be\label{4.6}
\|\nabla \bW\|_{\infty , B_{1/2} }^4 \le   c (1+ |\bLambda |^{12} ) \intl_{B_1} |\bW- \bW_{B_1} |^4 dx.
\ee
Applying Poincar\'e's inequality, we obtain 
\be\label{4.7}
\intmw_{B_\tau } |\bW- \bW_{B_\tau } |^4 dx \le c\tau^4 \|\nabla \bW\|^4_{\infty , B_{1/2 }} .
\ee
Combination  of  \eqref{4.6} and \eqref{4.7} yields the desired estimate.
\hfill \Beweisende

\vspace{0.3cm}
\hspace*{0.5cm}
In our discussion below we make use of the notion of the Morrey space. We say 
$\bF \in {\cal M}^{p, \lambda }_{\rm loc} (\Omega)$ if for all $K \subset \subset \Omega $
\[
[\bF]^p_{{\cal M}^{p, \lambda }, K} := \sup\bigg\{  r^{-\lambda } \intl_{B_r(x_0)} |\bF|^{p}  dx  \,\bigg|\, x_0\in K, 0<r\le \dist(K,\partial \Omega ) \bigg\}   <+\infty.  
\] 

\begin{lem}
Let $\bF \in {\cal M}^{2, \lambda }_{\rm loc} (\Omega)$ for some $1< \lambda \le 3$. 
For every  $0<\tau < \frac {1} {2}, 0<M < +\infty, K\subset \subset \Omega $ and $0<\alpha < 
\frac {\lambda -1} {2}$,  
there exist positive numbers 
 $\var _0= \var _0(\tau , M, K, \alpha )$ and $R _0 = R _0(\tau , M, K, \alpha ) < 
\dist(K,\partial \Omega ) $ 
such that, if  $\bB\in W^{ 1,\, 2}_{\rm loc} (\Omega )$ is a suitable weak solution to \eqref{3.0}, and 
for $x_0\in K$  and  $0<R\le R_0$ the condition    
\be\label{4.8}
|\bB_{R, x_0 }| \le M, \quad E (R, x_0) + R^{\alpha } \le \var _0
\ee 
is fulfilled, then 
\be\label{4.9}
E (\tau R, x_0) \le 2 \tau C_0 (1+M^{3}) (E (R, x_0) + R^{\alpha}), 
\ee
where $C_0>0$   stands for the constant which appears  on the right-hand side of \eqref{4.3}. 
\end{lem}

{\it Proof}  Assume the assertion of the Lemma is not true. Then there exist $0<\tau < \frac {1} {2}, 
0<M < +\infty , K \subset \subset \Omega $  and $0<\alpha < \frac {\lambda -1} {2}$  
together with a sequence 
$\bB^{(k)}  \in W^{ 1,\, 2}_{\rm loc} (\Omega )$ being suitable weak solutions to \eqref{3.0} 
as well as sequences  $x _k \in K, 
0<R_k < \dist(K,\partial \Omega ) $ and $\var _k \rightarrow 0$ as 
$k \rightarrow +\infty $ such that 
\be\label{4.10}
|\bB^{(k)} _{ R_k, x_k}| \le M, \quad E_k (R_k, x_k) + R_k^{\alpha } = \var _k
\ee 
and 
\be\label{4.11}
E_k (\tau R_k, x_k) > 2 \tau C_0 (1+M^{3}) (E_k (R_k, x_k)+ R_k^{\alpha } ).
\ee
Here we have used the notation
\[
E_k (r, x_0) = \bigg(\intmw_{B_r(x_0)} |\bB^{(k)} - \bB^{(k)}_{r, x_0}  |^4 dx\bigg)^{1/4},\quad 
x_0 \in K, 0<r\le \dist(K,\partial \Omega ).   
\]
Note that \eqref{4.10} yields $R_k \rightarrow 0$ as $k \rightarrow +\infty $. 

\hspace*{0.5cm}
Next, define
\begin{align*}
\bW_k (y) &= \frac {1} {\var _k} (\bB^{(k)} (x_k + R_k y) - \bB^{(k)} _{ R_k, x_k}),
\\     
\bF_k (y) &=  \bF(x_k + R_k y),\quad y \in B_1(0) 
\end{align*}
$( k\in \N)$.  Furthermore, set
\[
\mathscr{E}_k (\sigma ) = \bigg(\intmw_{B_\sigma } |\bW_k- (\bW_k)_{B_\sigma }|^4  dy\bigg)^{1/4},
\quad 0<\sigma \le 1.   
\]
Then \eqref{4.10} and \eqref{4.11}  turn into
\be\label{4.10-1}
|\bB^{(k)} _{ R_k, x_k}| \le M, \quad \mathscr{E}_k (1) + \frac { R_k ^{\alpha  }} {\var _k}= 1 ,
\ee 
and 
\be\label{4.11-1}
\mathscr{E}_k (\tau ) > 2 \tau C_0 (1+M^{3}) \Big(\mathscr{E}_k (1)+ 
\frac { R_k^{\alpha  }} {\var _k}\Big)= 2 \tau C_0 (1+M^{3})
\ee
respectively. 

\hspace*{0.5cm}
Using the chain rule,  we find that \eqref{3.0} transforms into 
\begin{align}
- \Delta \bW_k &= - 
\var _k \nabla \times ((\nabla \times  \bW_k) \times  \bW_k)
- \nabla \times ((\nabla \times  \bW_k) \times  \bB^{(k)}_{R_k, x_k}) 
\cr
& \qquad +\frac {R_k} {\var _k} \nabla \times  \bF_k \quad \mbox{in}\quad B_1.
\label{4.12}
\end{align}
Thus, $\bW_k \in W^{ 1,\, 2}(B_1)$ is a weak solution to \eqref{4.12}.  
Let $0<\sigma <1$. Using the transformation formula, noticing that $|\bB^{(k)}_{R_k, x_k}|\le M  $,  
the Caccioppoli-type inequality 
\eqref{3.4-1} with $r= R_k$ and $\rho = \sigma  R_k$  turns into 
\begin{align}
 \|\nabla \bW_k\|_{2, B_\sigma } &\le c(1-\sigma )^{-1} 
\Big((1+M) \mathscr{E}_k(1) + 
\var_k  \mathscr{E}_k(1)^2\Big)  + \frac {c R_k^{-1/2}} {\var _k  } \|\bF \|_{2,B_{R_k}(x_k)}.  
\label{4.12-1}
\end{align}
Observing \eqref{4.10-1}, and verifying 
\be\label{4.12-1a}
\frac {R_k^{-1/2} } {\var _k}\|\bF \|_{B_{R_k}(x_k)}  \le \frac { R_k^{(\lambda -1)/2}} {\var _k} 
[\bF]_{ {\cal M}^{2, \lambda }, K } 
\le  R_k^{(\lambda -1)/2-\alpha }   [\bF]_{ {\cal M}^{2, \lambda }, K }
\ee
from \eqref{4.12-1}, we get 
\be\label{4.12-2}
\|\nabla \bW_k\|_{2, B_\sigma } \le c(1-\sigma )^{-1} (M+ 1)  + 
c [\bF]_{ {\cal M}^{2, \lambda }, K }.   
\ee
In addition, in view of  \eqref{4.10-1} we estimate 
\be\label{4.12-3}
\|\bW_k\|_{4, B_1} = (\mes B_1)^{1/4}  \mathscr{E}_k (1) \le (\mes B_1)^{1/4}.  
\ee
From \eqref{4.12-2} and \eqref{4.12-3} it follows that $(\bW_k)$ is bounded in $W^{ 1,\, 2}(B_\sigma )$ for all $0<\sigma <1$ and bounded in $L^4(B_1)$.  
Thus, by means of reflexivity,  eventually passing to subsequences, 
we get $\bW\in  L^4(B_1)\cap  W^{ 1,\, 2}_{\rm loc} (B_1)$ and  $\bLambda \in \R^3$ such that 
\begin{align}
 \bW_k &\rightarrow \bW \quad \mbox{{\it weakly in}} \quad L^4(B_1) \quad \mbox{as} \quad k \rightarrow +\infty, 
\label{4.13}
\\
\bW_k &\rightarrow \bW \quad \mbox{{\it weakly in}} \quad W^{ 1,\, 2}(B_\sigma )
\quad \mbox{as} \quad k \rightarrow +\infty\quad \forall\, 0<\sigma <1, 
\label{4.13-1}
\\
\bB^{(k)} _{ R_k, x_k} &\rightarrow \bLambda  \quad\mbox{in}\quad \R^3 \quad 
\mbox{as}\quad k \rightarrow +\infty. 
\label{4.14}
\end{align}
On the other hand, by the compactness 
of the embedding $ W^{ 1,\, 2}(B_\sigma  ) \hookrightarrow L^4(B_\sigma )$ from \eqref{4.13} we infer 
\be\label{4.19}
\bW_k \rightarrow \bW \quad \mbox{{\it strongly in}} \quad L^{4}(B_\sigma ) \quad \mbox{as} \quad k \rightarrow +\infty \quad \forall\, 0<\sigma < 1.
\ee   
Accordingly,
\be\label{4.20}
\lim_{k \to \infty } \mathscr{E}_k (\sigma ) = \mathscr{E} (\sigma )\quad \forall\,  0<\sigma <1, \ee 
where $\mathscr{E}(\sigma )= {\D \bigg(\intmw_{B_\sigma } |\bW- \bW_{B_\sigma }  |^4 dy\bigg)^{1/4}}$. 
In particular, by the aid of \eqref{4.20} (with $\sigma =\tau $) from \eqref{4.11-1} we get 
\be\label{4.20-1}
\mathscr{E} (\tau ) \ge 2 \tau C_0 (1+ M^3). 
\ee 
In view of \eqref{4.12-1a} we have
 $$\frac {R_k} {\var _k} \|\bF_k\|_{2, B_{1}} = \frac {R_k^{-1/2} } {\var _k}\|\bF \|_{B_{R_k}(x_k)}\le R_k^{(\lambda -1)/2 -\alpha  } [\bF]_{{\cal M}^{2, \lambda }, K } \rightarrow 0 
 $$
 as $k \rightarrow +\infty $. Therefore, with the help of  \eqref{4.13}, \eqref{4.13-1}, \eqref{4.14} and \eqref{4.19},  letting $k \rightarrow+\infty  $ 
in \eqref{4.12}, we see that $\bW\in  W^{ 1,\, 2}_{\rm loc} (B_1)\cap  L^4(B_1)$ is a weak solution to 
\be\label{4.21}
-\Delta \bW = - \nabla \times ((\nabla \times  \bW)\times  \bLambda ) \quad \mbox{in} 
\quad B_1.
\ee  
As $|\bLambda |\le M$ appealing to  Lemma\,4.1, we find 
\be\label{4.22}
\mathscr{E}(\tau ) \le \tau C_0(1+ M^3) \mathscr{E}(1).    
\ee
On the other hand,  by virtue of   the lower semi continuity of the norm 
together with   \eqref{4.11-1} and \eqref{4.20} we get 
\begin{align*}
\mathscr{E}(1) &\le  \liminf_{k \to \infty } \Big(\mathscr{E}_k(1) + 
\frac {R_k^{\alpha } } {\var _k}\Big) \le    \frac {1} { 2\tau C_0(1+ M^3)} \lim_{k \to \infty } \mathscr{E}_k(\tau )
\\     
&= \frac {1} { 2\tau C_0(1+ M^3)} \mathscr{E}(\tau ).
\end{align*}
Estimating the right of \eqref{4.22} by the inequality we have just obtained, we see that  
$\mathscr{E}(\tau ) \le \frac {1} {2} \mathscr{E}(\tau )$ and hence  
$\mathscr{E}(\tau ) =0$, which contradicts to \eqref{4.20-1}. Whence, the assumption is not true, 
and this  completes   the proof of the Lemma. \hfill \Beweisende 

%
%
\section{Partial regularity} 
\setcounter{secnum}{\value {section}
\setcounter{equation}{0}
\renewcommand{\theequation}{\mbox{\arabic{secnum}.\arabic{equation}}}}

The aim of the present section is to prove the partial regularity of a suitable weak solution 
$\bB \in W^{ 1,\, 2}_{\rm loc} (\Omega )$ to system \eqref{3.0}, which will lead to the partial regularity of a suitable weak solution to the steady Hall-MHD system.  As we will see below,  the set $\Sigma  (\bB)$ of possible  singularities is given by means of
\begin{align}
\Sigma (\bB) &= \bigg\{ x_0\in \Omega \,\bigg|\,\liminf_{r \to 0^+} \frac {1} {r}\intl_{B_r(x_0)} |\nabla  \bB|^2 dx >0 \bigg\}  
\cr   
&\qquad\qquad \qquad\cup 
\Big\{x_0\in  \Omega   \,\Big|\,\sup_{r>0} |\bB_{r, x_0 }|=+\infty \Big\}.      
\label{5.0}
\end{align}

\begin{thm}
Let $\bF\in  {\cal M}^{2,\lambda }_{\rm loc} (\Omega ) $ for some $1<\lambda < 3$. 
Let $\bB\in  W^{ 1,\, 2}_{\rm loc} (\Omega )$ be a suitable weak solution to the system \eqref{3.0}.  Then, $\Omega \setminus \Sigma (\bB)$
 is an open set, on which $\bB$   is $\alpha$-H\"older continuous with repect to any  exponent $0<\alpha < \frac {\lambda -1} { 2}$.
\end{thm}

{\it Proof} Let $x_0 \in  \Omega \setminus  \Sigma (\bB)$ (cf. \eqref{5.0}). Set   $ d= \dist(x_0, \partial \Omega )$, and define $K= \overline{B_{d/2}  (x_0)} $. 
 Using Sobolev Poincar\'{e}'s inequality, we have
\[
\liminf_{r \to 0^+} \intmw_{B_r(x_0)} |\bB- \bB_{r, x_0 } |^4 dx \le 
c \liminf_{r \to 0^+} \bigg(\frac {1} {r}\intl_{B_r(x_0)} |\nabla  \bB|^2 dx\bigg)^2 =0  .   
\]
We set 
\[
M:= \sup_{0< r< d/2} |\bB_{r, x_0 }| + 3<+\infty. 
\]
Take $\tau  > 0$  such that 
\be\label{5.2}
2 \tau^{1-\alpha }   C_0 (1+ M^3) \le \frac {1} {2},\qquad \tau ^{\alpha } \le \frac {1} {2}. 
\ee 
Let $\var _0=\var_0(\tau , M, K, \alpha )$ and $ R_0=R_0(\tau , M, K, \alpha )$  denote the numbers according to Lemma\,4.2.  Then, we choose $0< R_1\le R_0$ such that 
\be\label{5.3}
E (R_1, x_0) + 2R^{\alpha }_1 <  \min\{\var _0, \var _1\},  
\ee
where $\var _1>0$ fulfills 
\be\label{5.3-1}
2\tau ^{-3} \var _1 \le 1. 
\ee
By the absolutely continuity of the Lebesgue integral there exists $ 0<\delta < \frac {d} {2}$ such that 
for all $y\in  \overline{B_\delta (x_0)}$
\begin{align}
& E(y, R_1) + 2R^{\alpha }_1  \le   \min\{\var _0, \var _1\},   
\label{5.4}
\\
& |\bB_{R_1, y}| \le \sup_{0<r< d/2} |\bB_{r, x_0} |+1 = M-2.
\label{5.5}
\end{align}
Fix $y\in \overline{B_\delta (x_0)}$.  We claim that for every $j \in  \N\cup \{0\}$  there holds 
\begin{align}
E(\tau ^{j} R_1, y) &\le 2^{- j} \tau^{\alpha j} E(R_1, y) + (1- 2^{-j} ) \tau ^{\alpha j} R_1^{\alpha },
\label{5.6}  
\\
|\bB_{ \tau ^j R_1, y}|&\le  M - 2^{-j+1}.
\label{5.7}
\end{align}
Clearly, for $j=0$, \eqref{5.6} is trivially fulfilled, while \eqref{5.7} holds in view of \eqref{5.5}.  

\hspace*{0.5cm}
Now, assume \eqref{5.6} and \eqref{5.7} are fulfilled for $j\in  N\cup \{0\}$.  Then \eqref{5.6} 
along with \eqref{5.3}  immediately implies
\be\label{5.8}
E(\tau ^{j} R_1, y) + \tau^{\alpha j } R_1^{\alpha }  \le \tau^{\alpha j} 
 (E(R_1, y) + 2R_1^\alpha)
\le \tau^{\alpha j} \min\{\var _0, \var _1\}.
\ee 
Now, both \eqref{5.7} and \eqref{5.8} imply  
\[
E(\tau ^{j} R_1, y) + \tau^{\alpha j} R_1^{\alpha}\le \var _0,\quad |\bB_{\tau ^j R_1.  y}|\le  M. 
\]
Thus, we are in a position to apply Lemma\,4.2 with $R= \tau ^j R_1$. This together with 
\eqref{5.2}  gives 
\begin{align}
E(\tau^{j+1} R_1, y) &\le 2 \tau C_0 (1+ M^3) (E(\tau^{j} R_1, y ) + 
\tau ^{\alpha j} R_1^{\alpha }) 
\cr
&\le \frac {1} {2} \tau^{\alpha} E(\tau^{j} R_1, y) + \frac {1} {2}\tau ^{\alpha (j+1)} R_1^{\alpha }
\cr
&\le 2^{-(j+1)} \tau^{\alpha (j+1) } E(R_1, y) + (1- 2^{-(j+1)} ) \tau ^{\alpha (j+1)} R_1^{\alpha }.
\label{5.9}
\end{align}
This proves \eqref{5.6} for $j +1$.  

\hspace{0.5cm}
It remains to show \eqref{5.7} for $j+1$.  First, from 
\eqref{5.6} along with  \eqref{5.4} we infer 
\be\label{5.10}
E(\tau ^{j} R_1, y) \le \tau^{\alpha j} (E(R_1, y) + R_1^{\alpha }) 
\le \tau ^{\alpha j} \var _1.  
\ee
Using triangle inequality and Jensen's inequality, we find 
\begin{align*}
  |\bB_{\tau^{j+1} R_1, y}| &\le |\bB_{ \tau^{j}  R_1, y}| + 
\Big|\bB_{\tau^{j+1}  R_1,  y}- \bB_{\tau^{j}  R_1, y}\Big|
\\     
&\le |\bB_{\tau^{j}  R_1, y}| + 2 \tau ^{-3} E( \tau ^{j} R_1, y ).  
\end{align*}
Estimating the first member on the right by using \eqref{5.7} and the second one 
by the aid of\eqref{5.10}  together with \eqref{5.3} and \eqref{5.3-1}, we obtain 
\begin{align*}
|\bB_{ \tau^{j+1}  R_1, y}|  & \le M- 2^{-j+1} + 2 \tau ^{-3} \tau ^{\alpha j} \var _1
\\     
&\le M-2^{-j+1} + 2^{-j}= M- 2^{-j}.  
\end{align*}
This completes the proof of \eqref{5.7} for $j+1$. Whence, the claim.

\hspace*{0.5cm}
From (\ref{5.6}) we get a constant $C_1>0$ such that 
\[
\bigg(\intmw_{B_r(y)} |\bB- \bB_{r, y} |^4  dx\bigg)^{1/4} \le  C_1 r^{\alpha } \quad \forall\, 0<r < \frac {d} {2}, \quad \forall\, y\in  \overline{B _\delta (x_0)}. 
\]  
Thus, by the well-known equivalence of the Campanato space and the H\"{o}lder space(see e.g.  \cite{cam} or Theorem 1.3 of \cite{gia}) we conclude 
\be\label{5.10-1}
\bB|_{\overline{B_\delta (x_0)}} \in  C^{\alpha }(\overline{B_\delta (x_0)}). 
\ee

\hspace{0.5cm}
Finally, we shall verify that $B_\delta (x_0)\subset  \Omega \setminus \Sigma (\bB)$.  Let $y\in  B_{\delta }(x_0) $ be arbitrarily chosen. 
Firstly, notice that $\sup_{0<r<d/2}  |\bB_{r, y}|<+\infty$ (see \eqref{5.10-1} ). Secondly, from the Cacciopploli-type inequality 
\eqref{3.4-1}  with $0<r < \frac {d} {2}$ and $\rho = \frac {r} {2}$ replacing $x_0$ by $y$ therein we deduce
\begin{align*}
&\bigg(\frac {1} {r} \intl_{B_{r/2}(y) }  |\nabla  \bB|^2  dx\bigg)^{1/2}  
\\
&\qquad \le 
c\Big\{(1+|\bB_{r, y} |) E(r,y) + E(r,y)^2\Big\} +  
c r^{(\lambda -1)/2} [\bF]_{{\cal M}^{2, \lambda }, K }.
\end{align*}
As $ \lim_{r \to 0^+} E(r,y) = 0$ and $\lambda >1$ the right-hand side of the above inequaliy tends to zero as  
$r \rightarrow 0^+$. Hence, $y\in  \Omega \setminus \Sigma (\bB)$. This completes the proof of the theorem. \hfill \Beweisende

\begin{rem}
By the result of  \cite[Chap. IV, 2.]{gia}  there holds 
\be\label{5.11}
{\cal H}^{\beta }(\Sigma (\bB))=0 \quad \forall\, \beta >1,
\ee
which implies that $\Sigma (\bB)$ has Hausdorff dimension  at most one. We don't know, however,  whether the one-dimensional Hausdorff measure of $\Sigma (\bB)$ is finite.  
\end{rem}

\vspace{0.3cm}
\hspace*{0.5cm}
As a consequence of Theorem\,5.1 we get the following partial regularity result for the steady Hall-MHD system.
\begin{thm}
Let $\bbf \in  L^{6/5}\cap L^2_{\rm loc} $  and $\bg \in L^2\cap  {\cal M}_{\rm loc} ^{2,\lambda } $
for some $1< \lambda \le 2$.  Let $(\bu, p, \bB)$ be a weak solution to the Hall-MHD system such that 
$\bB$ satisfies the local energy inequality \eqref{3.0-3} with $\bF = - \bB\times \bu + \bg$.  
Then,  for $\Sigma (\bB)$ defined by \eqref{5.0} we have that 
$\R^3 \setminus  \Sigma (\bB) $ is an open set such that  $\bB$ is $\alpha $-H\"older continuous on 
$\R^3 \setminus  \Sigma (\bB)$ for any $0<\alpha < \frac {\lambda -1} {2}$. 
\end{thm}

{\it Proof} Arguing similarly as in the proof of Theorem\,2.3, from $\bbf, \bg\in  L^2_{\rm loc} $  we deduce that $\bV = \bB + \bfomega \in \bW^{ 1,\, 2}_{\rm loc}$.  As $\bB \in W^{ 1,\, 2}$ we see that 
$\bfomega \in  W^{ 1,\, 2}_{\rm loc}$, and hence $\bu \in W^{ 2,\, 2}_{\rm loc}$.  By means of Sobolev's 
embedding theorem we find $\bu \in  L^\infty _{\rm loc} $.  This shows that 
$\bB\times \bu \in  L^6_{\rm loc} \subset  {\cal M}^{2,2}_{\rm loc}  $.  Consequently, 
$\bF = - \bB\times \bu + \bg \in  {\cal M}^{2,\lambda }_{\rm loc} $.  The assertion of the theorem 
is now an immediate consequence of Theorem\,5.1. \hfill \Beweisende

%
%
\section{Higher regularity} 
\setcounter{secnum}{\value {section}
\setcounter{equation}{0}
\renewcommand{\theequation}{\mbox{\arabic{secnum}.\arabic{equation}}}}

In Section\,5 we have proved the partial H\"older regularity of a suitable weak solution $(\bu, p, \bB)$ of the   Hall-MHD system for $\bbf $ and $\bg$ being sufficiently regular.  The aim of the present section is to show that if both $\bbf$ and $\bg$ are smooth, then $(\bu, p, \bB)$ is smooth in $\R^3\setminus  \Sigma (\bB)$.  
To prove this we first shall  establish a regularity result for the following linerized problem.  

\hspace*{0.5cm}
Let $\Omega \subset \R^3$ be an open set, and let $\bB\in C(\Omega )$. We consider the linear system
\be\label{6.1}
-\Delta \bA = -\nabla \times  ((\nabla \times \bA) \times \bB)+ \nabla \times \bF
 \quad \mbox{in}\quad \Omega. 
\ee
We have the following regularity result. 
\begin{thm}
Let $\bF\in  {\cal M}^{2, \lambda }_{\rm loc} (\Omega ) $ for some $1<\lambda <   3$. 
Let $\bA\in  W^{ 1,\, 2}_{\rm loc} (\Omega )$ be a weak solution to \eqref{6.1}. Then 
\be\label{6.2}
\bA \in  C^{\alpha }(\Omega )\quad \mbox{with} \quad \alpha =\frac {\lambda -1} {2}.  
\ee
\end{thm}
 
{\it Proof}  Let $x_0 \in  \Omega $. Set $d= \dist(x_0,\partial \Omega )$.  As $\bB$ is continuous, 
we get 
\[
M := \max_{y\in \overline{B_{d/2}(x_0)} } |\bB(y)| <+\infty.  
\]
Let $y \in  B_{d/4}(x_0)$ be fixed.   For the sake of notational simplicity in what follows we use the notation 
\[
S(R, y) = \bigg(\frac {1} {R} \intl_{B_R(y)} |\nabla \bA|^2  dx\bigg)^{1/2} ,\quad 0<R\le \frac {d} {4}. 
\]
Furthermore, by $\osc(f; x_0, R)$ we denote the oscillation of a continuous function  $f$  over  $B_R(x_0)$, which equals the supremum of $|f(x)-f(y)|$ taken over all $x,y\in  B_R(x_0)$.  As $\bB$ is continuous we may choose 
$0< R_0 < \frac {d} {4}$ such that 
\[
\sup_{B_{R} (y)} |\bB- \bB_{R, y} | \le  \osc (\bB; x_0, 2R)\le \frac {1} {2}\quad \forall\, 0<R\le R_0.  
\]

\hspace*{0.5cm}
Next, for  $0<R\le  R_0$  let  $\bZ \in  W^{ 1,\, 2}(B_R(y))$ denote the unique weak solution to 
\begin{align}
-\Delta \bZ &= -\nabla \times  ((\nabla \times \bZ) \times \bB_{R, y} )
\cr
&\qquad
-\nabla \times  ((\nabla \times \bA) \times (\bB- \bB_{R,y}) ) +
\nabla \times \bF  \quad \mbox{in}\quad B_{R}(y),
\label{6.4}
\\
\bZ &={\bf 0} \quad \mbox{on}\quad \partial  B_{R}(y).  
\label{6.5}
\end{align}
It can be easily checked that 
\begin{align}
R^{-1/2} \|\nabla \bZ\|_{2, B_R(y)}  
&\le 4\osc (\bB; x_0, 2R) R^{-1/2}  \| \nabla  \bA\|_{2, B_R(y)} 
+ 2 R^{-1/2} \|\bF\|_{2, B_R(y)}
\cr
&\le  4\osc (\bB; x_0, 2R) S(R,y) + 2R^{(\lambda-1)/2} [\bF]_{{\cal M}^{2, \lambda }, B_{d/2}(x_0)  }.  
\label{6.6}
\end{align}

Setting $\bW= \bA- \bZ$ and $\bLambda = \bA_{R, y} $, it follows that $\bW\in  W^{ 1,\, 2}(B_R(y))$ 
is a weak solution to \eqref{4.1} in $B_R(y)$.  Analogeously as Lemma\,4.1 one shows that  
get 
\be\label{6.7}
\bigg(\frac {1} {\tau R}\intl_{B_{\tau R}(y)}  |\nabla \bW|^2 dx\bigg)^{1/2}    
\le C_0 \tau (1+ M^2) \bigg(\frac {1} {R} \intl_{B_{ R}(y)}  |\nabla \bW |^2 dx\bigg)^{1/2} \ee for all $0< \tau  <1$. 

\hspace*{0.5cm}
Next, let $\frac {\lambda -1} {2} <\beta < 1$ be fixed.  Take $0<\tau < \frac{1}{2}$ such that 
\be\label{6.7-1}
C_0\tau ^{1-\beta } (1+M^2) \le \frac {1} {2}.   
\ee
Thus, using triangle inequality along with \eqref{6.6}, \eqref{6.7} and \eqref{6.7-1} ,
we get 
\begin{align*}
 S(\tau R, y) &\le 
\bigg( \frac{1}{\tau R}\intl_{B_{\tau R}(y)}  |\nabla  \bW _{\tau R, y} |^2 dx\bigg)^{1/2}        
+ \tau^{-1/2} R^{-1/2} \|\nabla \bZ\|_{2, R}  
\\     
& \le C_0 \tau (1+ M^2) S(R, y) + (1+ \tau^{-1/2}) R^{-1/2} \|\nabla \bZ\|_{2, R}   
\\
& \le \frac {1} {2}\tau^{\beta } S(R,y) + 4(1+ \tau^{-1/2})\osc (\bB; x_0, 2R) S(R,y) + C_1 R^{\alpha }, 
\end{align*}
where
\be\label{6.9}
C_1 = 2(1+ \tau ^{-1/2}) [\bF]_{{\cal M}^{2, \lambda }, B_{d/2}(x_0)},\quad \alpha = \frac {\lambda -1} {2}.
\ee
Take $0<R_1 \le R_0$ such that $4(1+ \tau^{-1/2})\osc (\bB; x_0, 2R_1) \le \frac {\tau ^\beta } {2} $. Then from the 
estimate above we deduce 
\be\label{6.10}
S(\tau R, y) \le \tau ^\beta S(R, y) + C_1 R^\alpha \quad \forall\, 0< R\le R_1. 
\ee
By using a routine iteration argument, we infer from \eqref{6.10} that 
that 
\begin{align*}
S(\tau ^k R_1, y) &\le \tau ^\beta S(\tau ^{k-1} R_1 ) + C_1 \tau^{\alpha k}  R_1^\alpha    
\\     
&\le \tau^{k\beta }  S(\tau ^{k-2} R_1 ) +   
C_1 \tau^{\alpha k}  (1+ \tau ^{\beta -\alpha }+ \ldots + \tau ^{(\beta -\alpha )(k-1)})  R_1^\alpha    
\\
& \le \tau ^{k \alpha } \Big(1+ C_1R_1^\alpha  \frac {1} {1- \tau ^{\beta -\alpha } }\Big).  
\end{align*}
Thus, there exists a constant $C_2>0$ such that 
\be\label{6.11}
S(R, y) \le C_2 R^\alpha \quad \forall\, 0< R \le \frac {d} {4}, \quad y \in  B_{d/4}(x_0).  
\ee
By the aid of the Poincar\'{e} inequality from \eqref{6.11} we obtain 
\be\label{6.13}
\bigg(\intmw_{B_R} |\bA- \bA_{R, y} |^2 dx \bigg)^{1/2}  \le c R^\alpha \quad 
\forall\, 0< R \le \frac {1} {4}, \quad y \in  B_{d/4}(x_0),
\ee
which leads to the desired H\"older regularity of $\bA$. \hfill \Beweisende

\vspace{0.3cm}
\hspace{0.5cm}
We are now in a position to prove the higher regularity for a continuous weak solution $(\bu, p, \bB)$ 
to the steady Hall-MHD system. More precisely, we have the following 
\begin{thm}
For $\bbf\in L^{6/5}$ and $\bg \in  L^2$ let $(\bu, p, \bB)$ be a weak solution to the steady 
Hall-MHD system.  Let $\Omega \subset \R^3$   be an open set such that $\bB$ is continuous in $\Omega $ and $\bbf, \bg\in  C^k(\Omega )$ \, $(k\in  \N\cup \{0\})$. Then $\bB, \bfomega \in C^{k, \alpha }(\Omega ) $ for all $0<\alpha < 1$. 
\end{thm}

{\it Proof}  First,  let us consider the case $k=0$.  As $\bbf, \bg \in L^\infty _{\rm loc}(\Omega ) $ 
by virtue of Theorem\,2.3 and Sobolev's embedding theorem we get 
\be\label{6.15-1}
\bV \in W^{ 1,\, q}_{\rm loc} (\Omega )\quad \forall\, 1\le q <+\infty,\quad 
\bV \in C^\alpha _{\rm loc}(\Omega )\quad \forall\, 0<\alpha <1.
\ee
With help of Sobolev's embedding theorem we see that 
$-\bB\times \bu + \bg\in {\cal M}^{2, \lambda } _{\rm loc} $ for all $0<\lambda <3$. Hence,  
from Theorem\,6.1 with $\bA=\bB$ we immediately get $\bB\in  C^{\alpha  }_{\rm loc}(\Omega )  $ 
for all $0<\alpha <1$.  As $\bfomega = \bV - \bB$ we infer  $\bfomega \in  C^\alpha _{\rm loc}(\Omega ) $ 
and since $\nabla \cdot \bu=0$ it follows $\bu \in C^{1, \alpha}_{\rm loc}(\Omega )$.  This completes the 
proof of the assertion in case $k=0$.
  
\hspace*{0.5cm}
Suppose $\bbf, \bg \in  C^k_{\rm loc} (\Omega )$ for some $k\in \N$.  From the proof above we immediately get 
$\bB, \bfomega\in  C^\alpha _{\rm loc} (\Omega ) $ for all $0<\alpha <1$. 
Now,  assume that  $\bB, \bfomega \in C^{j-1, \alpha }_{\rm loc}  (\Omega ) 
\cap W^{j,\, 2}_{\rm loc} (\Omega )$  for all $0<\alpha <1$  for some $j\in \{1, \ldots, k-1\}$.   
Let $\nu \in  \N^3$ be a multi-index with $|\nu |=j-1$.   Define $\bA= D^\nu \bB$ in $\Omega $. 
Applying $D^\nu $ to both sides of \eqref{1.3},  we  are  led to 
\begin{align}
 - \Delta \bA& = - \nabla \times ((\nabla \times \bA )\times \bB) 
\cr
 & \quad  - \suml_{ |\mu |\le j-2, \mu \le \nu } \nabla \times ((\nabla \times D^\mu \bB )\times D^{\nu -\mu } \bB) 
+\nabla \times  D^\nu (\bB \times \bu + \bg)
\cr
&= - \nabla \times ((\nabla \times \bA )\times \bB) + \nabla \times  \bG
\label{6.20}
\end{align}
in $\Omega $. By our assumption, we have $\bG \in  C^\alpha_{\rm loc}(\Omega ) $ for all $0<\alpha <1$.  
Applying the method of differences,  we  see that 
$\bA \in  W^{1+\theta ,\, 2}_{\rm loc} (\Omega )$ 
for every $0< \theta <1$.   Consequently, $\bB \in  W^{ j+\theta ,\, 2}_{\rm loc} (\Omega )$ 
for all $0<\theta <1$.  By virtue of Sobolev's embedding theorem it follows that 
\be\label{6.14}
\bB \in  W^{ j ,\, q} _{\rm loc}(\Omega ) \quad \forall\, 1\le q < 6.  
\ee
This shows that the $\nabla \times  \bG \in  L^2_{\rm loc}(\Omega ) $.  
Therefore,  
we are able to perform the method of difference  quotient  which yields 
$\bA \in  W^{ 2,\, 2}_{\rm loc} (\Omega ) $. Recalling our assumption,  having 
$ \bA\in  C^\alpha_{\rm loc}  (\Omega )$ for all $0<\alpha < 1$  by  the interpolation
inequality due to Kufner and Wannebo \cite{kuf}, 
we obtain $\bA\in  W^{ 1,\, q}_{\rm loc} (\Omega )$ for all $1\le q < +\infty $. 
This proves that 
\be\label{6.21}
\bB \in W^{ j+1, 2}_{\rm loc} (\Omega ) \cap    \bigcap_{1\le q<\infty } 
W^{ j,\, q}_{\rm loc} (\Omega ).     
\ee

Repeating the above argument with $|\nu |=j$ and $\bA= D^\nu \bB$,  we see that 
$\bA \in  W^{ 1,\, 2}_{\rm loc} (\Omega )$ is a weak solution to 
\be\label{6.22}
-\Delta \bA = -\nabla \times (\nabla \times  \bA)\times \bB) + \nabla \times \bG.
\ee
Thanks to \eqref{6.21} we have $\bG \in  {\cal M}^{2, \lambda } _{\rm loc}(\Omega ) $ for all $0<\lambda < 3$, so that Theorem\,6.1 yields that $\bA \in C^\alpha_{\rm loc}  (\Omega )$ for all 
$0<\alpha <1$, and  that  implies 
\be\label{6.23}
\bB \in  C^{ j,\, \alpha}_{\rm loc} (\Omega )\quad \forall\, 0<\alpha <1.     
\ee

\hspace*{0.5cm}
Furthermore, according to our  assumption  we have $\bV\in  W^{ j-1,\, q}_{\rm loc} (\Omega )$ for all $1\le q< +\infty $. 
Consequently, $- \bV \times \bu + \bbf + \bg \in  W^{ j-1,\, q}_{\rm loc} (\Omega )$ for all 
$1\le q <+\infty $. Hence using the Calder\'{o}n-Zygmund inequality,  from \eqref{1.6}  we deduce that  
$ \bV \in W^{ j,\, q}_{\rm loc}(\Omega)$ for all $1\le q <+\infty $.  
In particular, $- \bV \times \bu + \bbf + \bg \in  W^{ j,\, q}_{\rm loc} (\Omega )$ for all 
$1\le q < +\infty $.   Once, more employing Calder\'on-Zygmund's inequality,  we find 
$\bV \in  W^{j+1, q}_{\rm loc} (\Omega ) $ for all $1\le q < +\infty $, and with the help 
of Sobolev's embedding theoren we get $\bV\in  C^{j, \alpha }_{\rm loc} (\Omega ) $ 
for all $0<\alpha <1$. Finally, recalling 
$\bfomega= \bV- \bB $  in view of \eqref{6.21} and \eqref{6.23}, we conclude 
\be\label{6.24}
\bfomega  \in  W^{j+1, 2}_{\rm loc} (\Omega ) \cap    \bigcap_{0< \alpha <1} 
C^{ j,\, \alpha}_{\rm loc} (\Omega ).   
\ee

\hspace{0.5cm}
The desired regularity now follows from above by induction over $j= 0, \ldots, k$. 

\hfill \Beweisende

%
%
\section{Direct method and compacteness of the singular set} 
\setcounter{secnum}{\value {section}
\setcounter{equation}{0}
\renewcommand{\theequation}{\mbox{\arabic{secnum}.\arabic{equation}}}}

In this section we prove that a suitable weak solution becomes regular outside a sufficiently large ball, which is due to the 
decay property 
\be\label{7.1}
\lim_{R \to +\infty } \intl_{\{|x|>R\} } |\bB|^6 dx = 0.    
\ee

\hspace*{0.5cm}
First let us state an alternative Caccioppoli-type inequality for the system \eqref{3.0}.

\begin{lem}
Let $\bF\in  L^2$ and let $\bB\in \dot{W}^{ 1,\, 2}$  be a suitable  weak solution to the system 
\eqref{3.0} in $\R^3$.  
Then for every Ball $B_r= B_r(x_0) $ 
and $0< \rho < r$ there holds 
\begin{align}
& \frac {1} {r} \intl_{B_{\rho  } }|\nabla  \bB|^2  dx  
\cr
&\quad\le \frac {cr^2} {(r-\rho )^2}  \bigg\{1+\bigg(\intmw_{B_r} |\bB|^6 dx\bigg)^{1/3}\bigg\}    \bigg(\intmw_{B_r} |\bB- \bB_{r, x_0} |^3  dx\bigg)^{2/3} 
\cr 
&\qquad\qquad+ \frac {c} {r} \|\bF\|_{2, B_r}^2,    
\label{7.2}
\end{align}
where $c=\const>0$ denotes a universal 
constant.

\end{lem}

{\it Proof}  This can be easily achieved by estimating the integral $J $ on the right-hand side of  (\ref{3.1})  
by using H\"older's inequality and Young's inequality as follows
\begin{align*}
J &\le \frac {c} {(r-\rho )^2} \intl_{B_r} |\bB|^2 |\bB- \bB_{r, x_0} |^2 dx  + 
\frac {1} {2}\intl_{B_r} \zeta ^2 |\nabla \bB|^2 dx     
\\     
& \le  \frac {cr^3} {(r-\rho )^2} \bigg(\intmw_{B_r} |\bB|^6 dx\bigg)^{1/3}  
\bigg(\intmw_{B_r} \bB- \bB_{r, x_0} |^3 dx\bigg)^{2/3}   + 
\frac {1} {2}\intl_{B_r} \zeta ^2 |\nabla \bB|^2 dx.
\end{align*}
\hfill \Beweisende

\hspace*{0.5cm}
Using the well-known properties of harmonic functions,  one easily verifies the following 

\begin{lem}
Let $\bW \in  W^{ 1,\, 3/2}(B_{R/2}(x_0) )$ be harmonic  in $B_{R/2}(x_0) $. Then, there exists an absolute constant $C_0$ such that for all $ 0< \tau < \frac{1}{2}$
\be\label{7.0-1}
\bigg(\intmw_{B_{\tau R} } |\bW- \bW_{\tau R, x_0} |^3 dx\bigg)^{1/3} 
\le C_0\tau \bigg(\intmw_{B_{ R/2} } |\bW- \bW_{ R/2, x_0} |^3 dx\bigg)^{1/3}.   
\ee

\end{lem} 
\hfill \Beweisende

\hspace*{0.5cm}
In what follows, let $\bF\in  {\cal M}^{2,\lambda } $ for some $1<\lambda <3$,  i.\,e.
\[
[\bF]^2_{{\cal M}^{2,\lambda } } =\sup \bigg\{ r^{-\lambda } \intl_{B_r(x_0)} |\bF|^2 dx 
\,\bigg|\,x_0\in \R^3, 0<r\le 1\bigg\}<+\infty.   
\]
Clearly 
\be\label{7.0-2}
R^{-1/2}\|\bF\|_{2, B_R} \le \gamma _0 R^\alpha \quad \forall\, 0<R\le 1,   
\ee
where
\begin{align*}
\alpha = \frac {\lambda -1} {2},\quad \gamma _0 = [\bF]_{{\cal M}^{2, \lambda}}. 
\end{align*}
Furthermore, define
\begin{align*}
 S(r, x_0) &=  \bigg(\frac{1}{r}\intl_{B_r(x_0)} |\nabla \bB|^2 dx\bigg)^{1/2},  
\\     
E(r, x_0) &= \bigg(\intmw_{B_r(x_0)} | \bB- \bB_{r, x_0} |^3 dx\bigg)^{1/3},
\\
M(r, x_0) &= \bigg(\intmw_{B_r(x_0)} | \bB |^6 dx\bigg)^{1/6},\quad 0<r<+\infty, \,\, x_0\in \R^3.
\end{align*}
Fix $x_0\in \R^3$ and $0<R\le 1$. Assume that $M(R, x_0)\le 1$. Then,
from \eqref{7.2} with $r=R $ and $\rho =\frac {R} {2}$ along with \eqref{7.0-2} we deduce 
\be\label{7.2-1}
S(R/2, x_0) \le c (E(R, x_0 )+  \gamma _0 R^{\alpha }),
\ee
where $c>0$ denotes an absolute constant. 

\hspace*{0.5cm}
Let $\alpha  < \beta  < 1$ be fixed. Take $0<\tau < \frac {1} {4}$ such that 
\be\label{7.8-1}
2C_0 \tau ^{1-\beta }\le \frac {1} {2}, \quad \tau ^{\beta -\alpha }\le \frac {1} {2},\quad 
\tau^\alpha \le \frac {1} {2}.  
\ee

Let $\bZ \in  W^{ 1,\, 3/2}(B_{R/2} (x_0))$  denote a weak solution  to 
\begin{align}
- \Delta \bZ &= - \nabla \times ((\nabla \times \bB)\times \bB) + \nabla \times \bF 
\quad \mbox{in}\quad B_{R/2} (x_0),
\label{7.4}  
\\
\bZ &={\bf 0} \quad \mbox{on}\quad\partial B_{R/2} (x_0).
\label{7.5}
\end{align}
By the well-known $L^p$-theory of the Laplace equation we get 
\begin{align}
R^{-1} \|\nabla \bZ\|_{3/2, B_{R/2} (x_0)} &\le cR^{-1} \|(\nabla \times \bB)\times \bB \|_{3/2, B_{R/2}(x_0) }
+ c R^{-1}\|\bF\|_{3/2, B_R(x_0)}        
\cr
&\le c M(R, x_0) S(R/2, x_0) + c \gamma _0 R \alpha.
\label{7.6}
\end{align}
Esimating the left hand side from below  by using   Sobolev-Poincar\'{e}'s inequality and the right hand side 
from above by the aid of \eqref{7.2-1}, recalling that $M(R, x_0)\le 1$,   we are led to 
\begin{align}
\bigg(\intmw_{B_{R/2}(x_0) } |\bZ|^3 dx\bigg)^{1/3} 
\le C_1 (M(R, x_0) E(R, x_0) + \gamma _0 R^{\alpha }),       
\label{7.7}
\end{align}
where $C_1>1$ stands for an absolute constant. 

\hspace*{0.5cm}
Next, we ssume that 
\be\label{7.9-1}
3\tau ^{-1} C_1 M(R,x_0)\le \frac {1} {2}\tau^{\beta }.
\ee
We note here that \eqref{7.9-1} yields $M(R, x_0)\le 1$ and thus \eqref{7.2-1} remains true.
We make use of triangle inequality, then apply   \eqref{7.0-1} (note that $ \bW  $ is harmonic). This together with  \eqref{7.7} and  \eqref{7.8-1}  gives  
\begin{align*}
E(\tau R, x_0) &\le   \bigg(\intmw_{B_{\tau R}(x_0) } |\bW- \bW_{\tau R, x_0} |^3 dx\bigg)^{1/3}
+ 3\tau ^{-1} \bigg(\intmw_{B_{ R/2}(x_0) } |\bZ |^3 dx\bigg)^{1/3}
\\   
&\le 2 C_0 \tau   E(R, x_0) + 3\tau ^{-1} \bigg(\intmw_{B_{ R/2}(x_0) } |\bZ |^3 dx\bigg)^{1/3}
\\
&\le \frac {1} {2}\tau^{\beta }    E(R, x_0) + 3\tau ^{-1} C_1 M(R,x_0) E(R, x_0) + 
3\tau ^{-1} C_1  \gamma_0 R^{\alpha }, 
\end{align*} 
and observing \eqref{7.9-1},  we therefore obtain
\be\label{7.10}
E(\tau R, x_0) \le \tau^{\beta }    E(R, x_0) + C_2 \gamma _0 R^{\alpha },
\ee
where 
\[
C_2 = 3 \tau ^{-1}C_1. 
\]

\hspace{0.5cm}
Next, we shall estimate $|M(\tau R, x_0) | $.  By using triangle inequality and Sobolev-Poincar\'{e} inequality 
it follows that
\begin{align}
M(\tau R, x_0)  &\le |\bB_{\tau R, x_0}| +  c\tau^{-1/2} S(R/2, x_0)      
\cr     
&\le M(R, x_0)  + c\tau^{-1/2} S(R/2, x_0) + 2\tau^{-1} E(R,x_0)
\cr     
&\le M(R, x_0)  + C_3  (E(R,x_0) + \gamma _0 R^\alpha)
\label{7.10-1}
\end{align}
with a constant  $C_3>0$  depending on $\tau $ only. 

\hspace*{0.5cm}
Let $0< M_0 \le 1$ such that 
\be\label{7.15}
3\tau ^{-1} C_1 M_0=  \frac {1} {2}\tau^{\beta }.
\ee
Let $0<R_1\le 1$ chosen so that 
\be\label{7.16}
2C_3(2C_2 \tau ^{-\alpha } + 1) \gamma _0 R_1^\alpha \le \frac {M_0} {2}.
\ee
Since  $\bB\in L^6$, there exists $0< \rho _0 <+\infty $ such that 
\be\label{7.17}
M(R_1, x_0) +2C_3 E(R_1, x_0) \le (1+ 4C_3)M(R_1, x_0)\le \frac {M_0} {2}\quad\forall\, |x_0|\ge \rho _0. 
\ee
Let $x_0 \in \R^n, |x_0|\ge \rho _0 $. We claim that for every $k\in \N\cup \{0\}$
\begin{align}
E(\tau ^k R_1, x_0) &\le \tau^{\beta  k} E(R_1, x_0) + 
2(1- 2^{-k}) \tau ^{\alpha (k-1)}  C_2\gamma _0 R_1^\alpha,
\label{7.18}
\\
M(\tau ^k R_1, x_0) &\le  M(R_1, x_0) + 2(1- 2^{-k} ) 
\Big\{C_3 E(R_1, x_0)  + C_3(2C_2 \tau ^{-\alpha } + 1) \gamma _0 R_1^\alpha\Big\}. 
\label{7.19}
\end{align}

\hspace{0.5cm}
We prove the claim by using induction over $k\in \N\cup  \{0\} $.  

\hspace{0.5cm}
Firstly, note that for $k=0$  both 
\eqref{7.18} and \eqref{7.19} are trivially fulfilled. 
Assume \eqref{7.18} and \eqref{7.19} hold for $k \in \N \cup \{0\}$.  Observing 
\eqref{7.16} and \eqref{7.17}, the assumption \eqref{7.19} implies 
\be\label{7.20}
M(\tau ^k R_1, x_0) \le M_0.
\ee
By the choice of $M_0$ \eqref{7.20} yields \eqref{7.9-1} for $R= \tau ^k R_1$.  Hence from 
\eqref{7.10} with $R= \tau ^k R_1$ we infer 
\[
E(\tau ^{k+1}R_1, x_0 ) \le \tau ^\beta E(\tau ^k R_1, x_0) + C_2 \gamma _0 \tau ^{\alpha k} R_1^\alpha. 
\]
Now, estimating the first term by the assumption \eqref{7.18} taking into account \eqref{7.8-1},  we arrive at 
\begin{align*}
E(\tau ^{k+1}R_1, x_0 ) &\le 
\tau^{\beta  (k+1)} E(R_1, x_0) + 2 \tau ^{\beta -\alpha } (1- 2^{-k} ) \tau ^{\alpha k}   C_2\gamma _0 R_1^\alpha   + \tau ^{\alpha k} C_2 \gamma _0 R_1^\alpha
\\     
&\le \tau^{\beta  (k+1)} E(R_1, x_0)+ (2- 2^{-k})  \tau ^{\alpha k}   C_2\gamma _0 R_1^\alpha 
\end{align*}
which results in \eqref{7.18} for $k+1$.  

\hspace*{0.5cm}
It remains to verify \eqref{7.19} for $k+1$. In fact, by means of \eqref{7.10-1} with $R= \tau ^k R_1$  together with the assumption  
\eqref{7.18} and \eqref{7.19} we estimate 
\begin{align*}
& M(\tau^{k+1} R_1, x_0 ) 
\\
&\le M(\tau^{k} R_1, x_0 ) + C_3(E(\tau ^k R_1, x_0) + \gamma _0 \tau ^{\alpha k} 
R_1^{\alpha } )
\\     
&\le M(R_1, x_0 ) + 2(1- 2^{-k} ) 
\Big\{C_3 E(R_1, x_0)  + C_3(2C_2\tau ^{-\alpha } + 1) \gamma _0 R_1^\alpha\Big\}
\\
&\qquad +  C_3 ( \tau^{\beta k } E(R_1, x_0) + 2\tau ^{\alpha (k-1)} C_2\gamma _0 R_1^{\alpha }  + 
\gamma _0 \tau ^{\alpha k} R_1^{\alpha } )
\\     
&\le M(R_1, x_0 ) + 2(1- 2^{-k} ) 
\Big\{C_3 E(R_1, x_0)  + C_3(2C_2\tau ^{-\alpha } + 1) \gamma _0 R_1^\alpha\Big\}
\\
&\qquad +  2^{-k} (C_3 E(R_1, x_0) + C_3(2 C_2 \tau ^{-\alpha } +1)\gamma _0 R_1^{\alpha })
\\
& = M(R_1, x_0 ) + 2(1- 2^{-k-1} ) 
\Big\{C_3 E(R_1, x_0)  + C_3(2C_2\tau ^{-\alpha } + 1) \gamma _0 R_1^\alpha\Big\}.
\end{align*}
Whence, \eqref{7.19} for $k+1$. 

\vspace{0.3cm}
\hspace{0.5cm}
This implies  the following  
\begin{thm}
Let $\bF\in  {\cal M}^{2, \lambda}$, $ 1< \lambda < 3$.  Let $\bB \in  \dot{W}^{ 1,\, 2}$ be a suitable weak solution to \eqref{3.0}. Then there exists $\rho _0>0$ such that $\Sigma (\bB)\subset B_{\rho _0}$. 
\end{thm}

\hspace{0.5cm}
As a consequence of Theorem\,7.3 we get 

\begin{cor} Let $\bbf \in L^{6/5} \cap  L^2$ and $\bg \in L^2\cap  L^q$ for some $3< q< +\infty $. 
Let $(\bu, p, \bB)$ be a suitable weak solution to the steady Hall-MHD system. Then there exists 
$\rho _0 >0$ such that $\bB$ is H\"older contiuous in $\{x:|x|>\rho _0\}$. In particular, 
$\Sigma (\bB)$ is a compact set of Hausdorff dimension at most one. 

\end{cor}

{\it Proof} To prove the corollary  we only need to verify that $- \bB\times \bu + \bg\in  {\cal M}^{2,\lambda } $
for some $1< \lambda <3$. Then the assertion follows immediately from Theorem\,7.3 with 
$\bF=- \bB\times \bu + \bg\in  {\cal M}^{2,\lambda }$. 

\hspace*{0.5cm}
First,  using H\"older's inequality,  we find $\|\bg\|^2_{2, B_R} \le c R^{3(q-2)/q}\|\bg\|_q $ for every ball $B_R \subset \R^3$,    which implies $\bg \in {\cal M}^{2, 3(q-2)/q} $.  
Owing to $3<q<+\infty $ we have $1<\frac {3(q-2)} {q}< 3$. 
Next, as $\bV = \bB + \bfomega \in  L^6 + L^2$  and $\bu \in  L^6$,  we see that 
$-\bV\times \bu+ (\bbf+ \bg) \in  L^{3/2}  + L^{3}+ L^2$.  By Calder\'{o}n-Zygmund's inequality 
it follows that $\nabla \bV\in  W^{ 1,\, 3/2}+ W^{ 1,\, 3}+ W^{ 1,\, 2}$.  By means of Sobolev's embedding 
theorem we get 
\begin{align*}
\|\bfomega \|_{3, B_2(x_0)} &\le  \|\bV\|_{3, B_2(x_0)}+ \|\bB\|_{3, B_2(x_0)}
\\     
&\le c (\|\bfomega \|_2\|\bu\|_6 + \|\bB\|_6\|\bu\|_6 + \|\bbf+ \bg\|_2) 
\end{align*}
for all $x_0 \in  \R^3$, with an absolute constant $c>0$.  As $\nabla \cdot \bu=0$,  we obtain 
\[
\|\bu\|_{10, B_{1} (x_0)} \le  c (\|\bu\|_6(1+ \|\bfomega \|_2 +  \|\bB\|_6) + c\|\bbf+ \bg\|_2\quad
 \forall\, x_0 \in \R^3. 
\]
Accordingly, $\bB \times \bu\in {\cal M}^{2, 7/5} $. This completes   the proof of the corollary. \hfill \Beweisende

\vspace{0.5cm}
{\bf Acknowledgements}
Chae was partially supported by NRF grants 2006-0093854 and  2009-0083521, while Wolf has been supported by the Brain Pool Project of the Korea Federation of Science and Technology Societies  (141S-1-3-0022).


\begin{thebibliography}{99}
 \bibitem{ach} M. Acheritogaray, P. Degond, A. Frouvelle  and J-G. Liu, {\it Kinetic formulation and global existence for the Hall-Magnetohydrodynamic system, } Kinetic and Related Models, {\bf 4} (2011), pp. 901-918.
 
 \bibitem{caf}L. Caffarelli, R. Kohn and L. Nirenberg, {\it Partial Regularity of suitable weak solutions of the Navier-Stokes equations,} Comm. Pure Appl. Math. {\bf 35} (1982), pp. 771-831.
 \bibitem{cam}S. Campanato, {\it Equazioni ellittiche del secondo ordine e spazi $ \mathcal{L}^{2, \lambda}$}, Ann. Mat. Pura e Appl. {\bf 69} (1965), pp. 321-380.
 
   \bibitem{cha1}D. Chae, P. Degond and J-G. Liu, {\it Well-posedness for Hall-magnetohydrodynamics}, Ann. Inst. Henri Poincare-Analyse Nonlineaire {\bf 31} (2014),  pp. 555-565.
   
     \bibitem{cha2}
D. Chae and J. Lee, {\it On the blow-up criterion and small data global existence for the Hall-magnetohydrodynamics.}  J. Differential Equations  {\bf 256}  (2014), no. 11, pp. 3835-3858.
 
\bibitem{cha3}
D. Chae and M. Schonbek, {\it On the temporal decay for the Hall-magnetohydrodynamic equations.} J. Differential Equations {\bf 255} (2013), no. 11, pp. 3971-3982.

\bibitem{cha4}
D. Chae and S. Weng, {\it Singularity formation for the incompressible Hall-MHD equations without resistivity,} Ann. Inst. Henri Poincare-Analyse Nonlineaire, (in print) http://dx.doi.org/10.1016/j.anihpc.2015.03.002.

  \bibitem{dre}
J. Dreher, V. Runban and R. Grauer, {\it Axisymmetric flows in Hall-MHD: a tendency towards finite-time singularity formation.} Physica Scripta  {\bf 72} (2005), pp. 451-455.  

\bibitem{fan}
J. Fan,  S. Huang and G. Nakamura, {\it  Well-posedness for the axisymmetric incompressible viscous Hall-magnetohydrodynamic equations.}  Appl. Math. Lett. {\bf 26} (2013), no. 9, pp. 963-967.

 \bibitem{for}T.G. Forbes, {\it Magnetic reconnection in solar flares,}  Geophys. Astrophys. Fluid Dyn.
 {\bf 62} (1991), pp.15-36.
 
\bibitem{gia}
M.\,Giaquinta, {\it  Multiple integrals in the calculus of variations and nonlinear elliptic systems}, Ann. of Math. studies {\bf 105}, Princeton Univ. press, Princeton, New Jersey,  1983.

\bibitem{hom}
H. Homann and  R. Grauer, {\it Bifurcation analysis of magnetic reconnection in Hall-MHD
systems.} Physica D {\bf  208} (2005), pp. 59-72.

\bibitem{kuf}
A. Kufner and A. Wannebo, {\it An interpolation inequality involving H\"{o}lder norms}, 
Georgian Math. J. {\bf 2}, no. 6 (1995), pp. 603-612.

  \bibitem{lad} O. A. Ladyzehnskaya and G. A. Seregin, {\it On partial regularity of sutiable weak solutions to the three-dimensional Navier-stokes equations,} J. Math. Fluid Mech. {\bf 1}, no. 4 (1999), pp. 356-387.
   
    \bibitem{lig}
M.\,J. Lighthill, {\it Studies on magneto-hydrodynamic waves and other anisotropic wave motions}, Philos. Trans. R. Soc. Lond. Ser. {\bf A 252}  (1960), pp. 397-430. 

\bibitem{lin}F. H. Lin, {\it A new proof of the Caffarelli-Kohn-Nirenberg theorem,} 
Comm. Pure Appl. Math. {\bf 51}  (1998), pp. 241-257.

\bibitem{pol} J. M. Polygiannakis and X. Moussas, 
{\it A review of magneto-vorticity induction in Hall-MHD plasmas,}
Plasma Phys. Control \& Fusion {\bf 43} (2001), pp.195-221.

\bibitem{sch}V. Scheffer, {\it Partial regularity of solutions to the Navier-Stokes equations,} Pacific J. Math. 
{\bf 66} (1976), pp. 535-552.

\bibitem{sha} D.A. Shalybkov, V.A. Urpin, {\it The Hall effect and the decay of magnetic fields,}  
Astron. Astrophys. (1997), pp. 685-690.

 \bibitem{miu}H. Miura and D. Hori, {\it Hall effects on local structure in decaying MHD turbulence,} J. Plasma Fusion Res. {\bf 8} (2009), pp. 73-76.
 
    \bibitem{sim}A. N. Simakov and L. Chac\'{o}n, {\it Quantitative, analytical model for magnetic reconnection in Hall magnetohydrodynamics,} Phys. Rev. Lett. {\bf 101} (2008), 105003.
 
 \bibitem{war}M. Wardle, {\it Star formation and the Hall effect, } Astrophys. Space Sci. {\bf 292},
   (2004) pp. 317-323. 
 
 \bibitem{wol}J. Wolf, {\it On the local regularity of suitable weak solution to the generalized Navier-Stokes equations,}  Ann. Univ. Ferrara, DOI 10.1007/s11565-014-0203-6.
 

\end{thebibliography}
\end{document}